\documentclass[11pt]{article}
                         \usepackage[latin1]{inputenc}
                         \usepackage{amsmath,amsthm,amssymb}
                         \newtheorem{teo}{Theorem}

                         \newtheorem{lemma}{Lemma}

                         \def\endproof{\nolinebreak\hfill\rule{2mm}{2mm}}
                         
                         \def\eq#1{(\ref{#1})}
                         
                         \def\neweq#1{\begin{equation}\label{#1}}
                         \def\endeq{\end{equation}}

                         \def\phi{\varphi}
                         \def\RR{{\mathbb R} }

                         \def\di{\displaystyle}
                         \def\ri{\rightarrow}
                         
                         \title{\sc Nonradial blow-up solutions of sublinear
                         elliptic equations with
                         gradient term}
                         \author{Marius GHERGU and
                         Vicen\c tiu R\u ADULESCU\thanks{Correspondence
                         address: Vicen\c tiu R\u adulescu,
                         Department of Mathematics, University of Craiova,
                         1100 Craiova, Romania, fax: +40-251.41.16.88. E-mail:
                         {\tt radulescu@inf.ucv.ro}}\\
                          \small Departament of Mathematics,
                                 University of Craiova,
                                         1100 Craiova, Romania}
                         \date{}
                         
\begin{document}
                         \baselineskip12pt

                         \maketitle
                         \noindent{\small {\bf Abstract}.
                         Let $f$ be a continuous and
                         non-decreasing function such that $f>0$ on $(0,\infty
                         )$,
                          $f(0)=0$, $\sup _{s\geq 1}
                         f(s)/s< \infty$ and let $p$ be a non-negative
                         continuous
                         function. We study the existence and nonexistence of
                         explosive solutions to the equation
                         $\Delta u+|\nabla u|=p(x)f(u)$ in $\Omega,$ where
                         $\Omega$ is either a
                         smooth bounded domain or $\Omega=\RR^N$. If $\Omega$
			 is bounded we prove that the above problem has never
			 a blow-up boundary solution. Since $f$
                         does not satisfy the Keller-Osserman growth condition at infinity,
                         we supply in the case $\Omega=\RR^N$ a necessary and
                         sufficient condition
                         for the existence of a positive solution that blows up
                         at infinity.
                         \medskip


                         }

                         \section{Introduction and the main results}
                         Explosive solutions of semilinear elliptic equations
                         have been studied
                         intensively in the last few decades.
                         Most of such studies have been concerned with
                         equations of the type
                         $$\Delta u=g(x,u),$$
                         in which the function $g$ takes various forms
                         (see \cite{bm, cia, cr1, cr2, cr3, gal, mar}
                         and their references).

                         In this paper we study an elliptic problem involving
                         a sublinear nonlinearity.
                         Due to the lack of the Keller-Osserman condition
                         \cite{k, o}, we find
                         a necessary and sufficient condition satisfied by the
                         potential so that our
                         problem admits a nonradial
                         solution blowing up at infinity. More precisely, we
                         consider the equation
                         \neweq{unu}
                          \left\{\begin{tabular}{ll}
                         $\Delta u+|\nabla u|=p(x)f(u)$ \qquad & $\mbox{ in }\
                         \Omega,$\\
                         $u\geq 0$ \qquad & $\mbox{ in }\ \Omega,$\\
                         \end{tabular} \right.
                         \endeq
                         where
                         $\Omega\subset\RR^N$ ($N\geq 3$) is either a smooth
                         bounded domain or the whole space.

			 The presence of the gradient term can have significant influence
			 on the existence of a solution, as well as on its asymptotic behavior.
			 Problems of this type appear in stochastic control theory and have been
			 first studied by Lasry and Lions \cite{la}. The corresponding
			 parabolic equation was considered in Quittner \cite{qu}. We also refer to
			 Bandle and Giarrusso \cite{bg, gia}
                         who established existence results and the asymptotic
                         behavior of solutions
                         for semilinear elliptic equations in bounded domains
                         containing gradient term (see also \cite{aw1} for
                         another class of
                         nonlinear elliptic problems involving gradient term).

                         Throughout this paper we assume that $p$ is a
                         non-negative function such that
                         $p\in C^{0,\alpha}(\overline\Omega )\,(0<\alpha<1)$ if
                         $\Omega$ is bounded, and
                         $p\in C^{0,\alpha}_{\rm loc}(\RR^N)$, otherwise. The
                         non-decreasing
                         non-linearity $f$ fulfills

			 \smallskip
			\noindent $\di (f1)\qquad f\in C^{0,\alpha}_{\rm
                         loc}[0,\infty),\;\;f(0)=0\,$ and
                         $\,f>0\,$ on $\,(0,\infty )$.

			 \smallskip
			 We also assume that $f$ is sublinear at infinity, in the sense that

                         \smallskip
                         \noindent $\di (f2)\qquad
                         \;\di\Lambda\equiv \sup _{s\geq 1}\frac{f(s)}{s}<
                         \infty.$

                         \smallskip
                         Cf. V\'eron \cite{ver}, the non-decreasing
                         non-linearity $f$ is
                         called an absorption term.

                         A solution $u$ of the problem \eq{unu} with
                         $u(x)\to\infty$ as $ {\rm dist}\,
                         (x,\partial\Omega)\to 0 $ (if
                         $\Omega$ is bounded) is called a {\it large
                         (explosive, blow-up)}
                         solution. If $\Omega = \RR^N,$ this condition can be
                         rewritten as
                         $u(x)\to\infty$ as $|x|\to\infty$. In this latter case
                         such solution
                         is called an {\it entire large (explosive)} solution. In terms of the dynamic
			 programming approach, an explosive solution of \eq{unu} corresponds to
			 a value function (or Bellman function) associated to an infinite
			 exit cost (see \cite {la}).

                         We note that in \cite{gr} it is studied the existence
                         and nonexistence
                         of large solutions for the corresponding
                         system to \eq{unu} where the coefficients are radial
                         functions.

			 \medskip If
                         $\Omega$ is bounded we prove the following non-existence result.
                         \medskip

                         \begin{teo}\label{th1}
                         Suppose $\Omega\subset\RR^N$ is a smooth bounded
                         domain.
                         Then \eq{unu} has no positive large solution in
                         $\Omega$.
                         \end{teo}

                         \medskip

                          Following Bandle and Giarrusso \cite{bg}, in
                         the presence of  the
                         Keller-Ossermann condition on $f,$ equation \eq{unu}
                         may have positive large solutions.
                         \medskip

                         Next, we consider problem \eq{unu}
                         when $\Omega=\RR^{N}$.
                         For all $r\geq 0$ we set
                         $$\phi(r)=\max \limits _{|x|=r}p(x),\qquad\psi(r)=\min
                         \limits
                         _{|x|=r}p(x),\qquad
                         \mbox{and}\qquad h(r)=\phi(r)-\psi(r).$$
                         We suppose that
                         \neweq{doi}
                         \int\limits_{0}^{\infty}rh(r)\Psi(r)dr<\infty,
                         \endeq
                         where
                         $$\Psi(r)=\exp\left(\Lambda_N\int\limits_{0}^{r}s\psi(s)ds\right),
                         \qquad  \Lambda_N=\frac{\Lambda}{N-2}.\\$$
                          Obviously, if $p$ is radial then $h\equiv 0$ and
                         \eq{doi} occurs.
                         Assumption \eq{doi} shows that the variable potential
                         $p(x)$ has a slow variation.
                         An example of nonradial potential for which \eq{doi}
                         holds is
                         $\di p(x)=\frac{1+|x_1|^2}{(1+|x_1|^2)(1+|x|^2)+1}.$
                         In this case
                         $\di\phi(r)=\frac{r^2+1}{(r^2+1)^2+1}$ and
                         $\di\psi(r)=\frac{1}{r^2+2}.$
                         If $\Lambda_N=1,$ by direct computation we get
                         $rh(r)\Psi(r)=O\left(r^{-2}\right)$ as $r\ri\infty$ and
                         so \eq{doi} holds.
                         \medskip

                          Our analysis will be developed under the
                         basic assumption \eq{doi}.
                         \medskip

                         \begin{teo}\label{th2} Assume $\Omega =\RR^N$ and
                         $p$ satisfies \eq{doi}.
                         Then \eq{unu} has positive entire large solution if
                         and only if
                         \neweq{apatru}
                         \int\limits_{1}^{\infty}e^{-t}t^{1-N}
                         \int\limits_{0}^{t}e^ss^{N-1}\psi(s)dsdt=\infty.
                         \endeq
                         \end{teo}
                         \medskip

                         \noindent\textbf{Remark 1.} Since
                         $\,\di
                         \int_0^re^rt^kdt=k!e^r\left(\sum_{s=1}^{k}(-1)^{k-s}\frac{t^s}{s!}\right)
                         -(-1)^kk!,$
                         for all integers $k\geq 1,$
                         we can give some examples of potentials $p$ that
                         verify both conditions \eq{doi}
                         and \eq{apatru}. In the case where $\Lambda_N=1$ such
                         functions are
                         \medskip

                         \noindent(i) $\di p(x)=1+|x|^m+|x_1|e^{-|x|^{m+2}},$
                         $m>0.$
                         \smallskip

                         \noindent(ii) $\di
                         p(x)=\frac{1+|x_1|g(|x|)e^{-|x|}}{1+|x|},$
                         $\,g\in C^{0,\alpha}_{\rm loc}[0,\infty)\cap
                         L^1[0,\infty),\,g\geq 0.$
                         \medskip

                         \noindent\textbf{Remark 2.}
                         We point out that a solution
                         of \eq{unu} may exist even if condition \eq{doi}
                         fails, as shown in what follows.
                          Define
                         $$\di
                         p(x)=2|x|^2+6x_1^2+\sqrt{|x|^2+3x_1^2}+N+1,\;x\in\RR^N.$$
                         and $f(t)=2t.$
                         For this choice of $p$ and $f,$ the equation \eq{unu}
                         has the nonradial entire
                         large solution
                         $u(x)=e^{|x|^2+x_1^2}$. In this case $h(r)=6r^2+r,$ so
                         \eq{doi} fails to hold.

			 \smallskip
			 The above results also apply to problems on Riemannian manifolds
			if $\Delta$ is replaced by the Laplace--Beltrami operator
			$$\Delta_B=\frac{1}{\sqrt c}\,\frac{\partial}{\partial x_i}\left(\sqrt
			c\, a_{ij}(x)
			\frac{\partial}{\partial x_i}\right)\, ,\qquad c:=\mbox{det}\,
			(a_{ij})\,,$$
			with respect to the metric $ds^2=c_{ij}\,dx_idx_j$, where $(c_{ij})$ is
			the inverse of $(a_{ij})$. In this case our results apply to concrete
			problems arising
			in Riemannian geometry. For instance, (cf.
			 Loewner-Nirenberg \cite{ln}) if $\Omega$ is replaced by the standard
			$N$--sphere $(S^N,g_0)$, $\Delta$ is the Laplace-Beltrami operator
			$\Delta_{g_0}$,
			$a=N(N-2)/4$, and $f(u)=(N-2)/[4(N-1)]\, u^{(N+2)/(N-2)}$, we find the
			prescribing scalar curvature equation with gradient term.

                         \section{Proofs}
                         \subsection{Proof of Theorem 1}
                         Suppose by contradiction that \eq{unu} has a positive
                         large solution $u$ and define
                         $v(x)=\ln(1+u(x)),\;x\in\Omega.$ It follows that $v$
                         is positive
                         and $v(x)\ri\infty$ as $\,{\rm dist}\,
                         (x,\partial\Omega)\to 0$.
                         We have
                         $$\di \Delta v=\frac{1}{1+u}\Delta
                         u-\frac{1}{(1+u)^2}|\nabla
                         u|^2\qquad\mbox{in}\;\Omega$$
                         and so
			 $$\di\Delta v\leq
                         p(x)\frac{f(u)}{1+u}\leq\|p\,\|_{\infty}\frac{f(u)}{1+u}\leq
                         A\qquad
                         \mbox{ in }\Omega,$$
                         for some  constant $A>0$.
                         Therefore
                         $$\di \Delta(v(x)-A|x|^2)<0,\qquad\mbox{for
                         all}\;\:x\in\Omega.$$
                         Let $w(x)=v(x)-A|x|^2,\;x\in\Omega$. Then $\di\Delta
                         w<0$ in $\Omega.$ Moreover,
                         since $\Omega$ is bounded, it follows that
                         $w(x)\ri\infty$ as ${\rm dist}(x,\partial\Omega)\ri
                         0.$

                         Let $M>0$ be arbitrary. We claim that $w\geq M$ in
                         $\Omega$.
                         For all $\delta>0,$ we set
                         $$\di \Omega_{\delta}=\{x\in\Omega\,;\;{\rm
                         dist}(x,\partial\Omega)>\delta\}.$$
                         Since $w(x)\ri\infty$ as ${\rm
                         dist}(x,\partial\Omega)\ri 0,$ we
                         can choose $\delta>0$ such that
                         \begin{equation}\label{awm}\di w(x)\geq M \qquad\mbox{
                         for all }\;
                         x\in\Omega\setminus\Omega_{\delta}.
                         \end{equation}
                         On the other hand,
                         $$\begin{tabular}{ll}
                         $\di -\Delta (w(x)-M)>0$ \qquad & ${\rm in}\
                         \Omega_{\delta},$\\
                         $\di\qquad\qquad w(x)-M\geq 0$ \qquad & ${\rm on}\
                         \partial\Omega_{\delta}.$\\
                         \end{tabular}$$
                         By the maximum principle we get
                         $w(x)-M\geq 0$ in $\Omega_{\delta}$. So, by
                          \eq{awm},
                          $w\geq M$ in $\Omega.$
                         Since $M>0$ is arbitrary, it follows that
                         $w\geq n$ in $\Omega$,
                         for all $n\geq 1$. Obviously, this is a
                         contradiction and the proof is now complete.
                         \endproof

                         \medskip
                         \subsection{Proof of Theorem 2}
                         Several times in the proof of Theorem 2 we shall apply
                         the  following inequality
                         \neweq{doispe}
                         \int\limits_{0}^{r}e^{-t}t^{1-N}
                         \int\limits_{0}^{t}e^ss^{N-1}g(s)dsdt\leq
                         \frac{1}{N-2}
                         \int\limits_{0}^{r}tg(t)dt, \qquad\forall\;r>0,
                         \endeq
			  for any continuous function
                          $g:[0,\infty)\ri [0,\infty)$.
                         Indeed, using an integration by parts in the left hand side
    			 we obtain
                         $$\begin{tabular}{ll}
                         $\di\int\limits_{0}^{r}e^{-t}t^{1-N}
                         \int\limits_{0}^{t}e^ss^{N-1}g(s)dsdt$ &
                         $\di\leq\int\limits_{0}^{r}t^{1-N}
                         \int\limits_{0}^{t}s^{N-1}g(s)dsdt$\\
                         &$\di=
                         \frac{1}{2-N}\int\limits_{0}^{r}\left(t^{2-N}\right)'
                         \int\limits_{0}^{t}s^{N-1}g(s)dsdt$\\
                         &
                         $\di=\frac{1}{2-N}\,r^{2-N}\int\limits_{0}^{r}t^{N-1}g(t)dt+
                         \frac{1}{N-2}\int\limits_{0}^{r}tg(t)dt$\\
                         &$\di\leq \frac{1}{N-2}\int\limits_{0}^{r}tg(t)dt,$
                         \end{tabular}$$
                         so \eq{doispe} follows.
                         \medskip

                         \noindent{\sc Necessary condition}.
                         Suppose that \eq{doi} fails and the equation \eq{unu}
                         has a positive entire
                         large solution $u$. We claim that
                         \neweq{phi}
                         \di\int\limits_{1}^{\infty}e^{-t}t^{1-N}
                         \int\limits_{0}^{t}e^ss^{N-1}\phi(s)dsdt<\infty.
                         \endeq

                         \noindent We first recall that $\phi=h+\psi.$ Thus
			 $$\begin{tabular}{ll}
                         $\di\int\limits_{1}^{\infty}e^{-t}t^{1-N}
                         \int\limits_{0}^{t}e^ss^{N-1}\phi(s)dsdt$&$\di=\int\limits_{1}^{\infty}e^{-t}t^{1-N}
                         \int\limits_{0}^{t}e^ss^{N-1}\psi(s)dsdt$\\
                         &$\di\;\;\;+\int\limits_{1}^{\infty}e^{-t}t^{1-N}
                         \int\limits_{0}^{t}e^ss^{N-1}h(s)dsdt.$\\
                         \end{tabular}$$
                         By virtue of \eq{doispe} we find
                         $$\begin{tabular}{ll}
                         $\di\int\limits_{1}^{\infty}e^{-t}t^{1-N}
                         \int\limits_{0}^{t}e^ss^{N-1}\phi(s)dsdt$&$\di\leq\int\limits_{1}^{\infty}e^{-t}t^{1-N}
                         \int\limits_{0}^{t}e^ss^{N-1}\psi(s)dsdt+\frac{1}{N-2}\int\limits_{0}^{\infty}th(t)dt$\\
                         &$\di\leq\int\limits_{1}^{\infty}e^{-t}t^{1-N}
                         \int\limits_{0}^{t}e^ss^{N-1}\psi(s)dsdt$\\
                         &$\di\;\;\;+\frac{1}{N-2}\int\limits_{0}^{\infty}th(t)\Psi(t)dt.$\\
                         \end{tabular}$$
                         Since $\di\int\limits_{1}^{\infty}e^{-t}t^{1-N}
                         \int\limits_{0}^{t}e^ss^{N-1}\psi(s)dsdt<\infty,$ by
                         \eq{doi}
                         we deduce that \eq{phi} follows.

                         Now, let $\bar u$ be the spherical average of $u,$
                         i.e.,
                         $$\bar
                         u(r)=\frac{1}{\omega_Nr^{N-1}}\int\limits_{|x|=r}
                         u(x)d\sigma_x,\qquad r\geq 0,$$
                         where $\omega_N$ is the surface area of the unit
                         sphere in $\RR^N$. Since $u$
                         is a positive entire large solution of \eq{unu} it
                         follows that $\bar u$
                         is positive and $\bar u(r)\ri\infty$ as $r\ri\infty.$
                         With the change of variable $\,x\ri ry,$ we have
                         $$\di\bar u(r)=\frac{1}{\omega_N}\int\limits_{|y|=1}
                         u(ry)\,d\sigma_y,\qquad r\geq 0$$
                         and
                         \begin{equation}\label{acinci}
                         \di\bar u'(r)=\frac{1}{\omega_N}\int\limits_{|y|=1}
                         \nabla u(ry)\cdot y\,d\sigma_y,\qquad r\geq 0.
                         \end{equation}
                         Hence
                         $$\di \bar u'(r)=
                         \frac{1}{\omega_N}\int\limits_{|y|=1}\frac{\partial
                         u}{\partial r}
                         (ry)\,d\sigma_y=
                         \frac{1}{\omega_Nr^{N-1}}\int\limits_{|x|=r}\frac{\partial
                         u}{\partial r}
                         (x)\,d\sigma_x,$$
                         that is
                         \begin{equation}\label{asase}
                         \bar
                         u'(r)=\frac{1}{\omega_Nr^{N-1}}\int\limits_{B(0,R)}
                         \Delta u(x)\,dx,\qquad \mbox{for all}\;\;r\geq 0.
                         \end{equation}

                         Due to the gradient term $|\nabla u|$ in \eq{unu}, we
                         cannot infer that
                         $\Delta u\geq 0$ in $\RR^N$ and so we cannot expect that
                         $\bar u'\geq 0$
                         in $[0,\infty)$. We define the auxiliary function
                         \begin{equation}\label{aopt}
                         \di U(r)=\max_{0\leq t\leq r}\bar u(t),\qquad r\geq 0.
                         \end{equation}
                         Then $U$ is positive and
                         non-decreasing.
                         Moreover, $U\geq\bar u$ and $U(r)\ri\infty$ as
                         $r\ri\infty$.

                         The assumptions $(f1)$ and $(f2)$ yield
                         $\di f(t)\leq \Lambda(1+t),$ for all $t\geq 0.$ So, by
                          \eq{acinci}
                         and \eq{asase},
                         $$\begin{tabular}{ll}
                         $\di \bar u''+\frac{N-1}{r}\,\bar u'+\bar
                         u'\!\!\!$&$\di\leq\;\frac{1}{\omega_Nr^{N-1}}
                         \int\limits_{|x|=r}\left[\Delta u(x)+|\nabla
                         u|(x)\right]d\sigma_x$\\
                         &$=\;\di \frac{1}{\omega_Nr^{N-1}}
                         \int\limits_{|x|=r}p(r)f(u(x))d\sigma_x$\\
                         &$\di\leq \;\Lambda \phi(r)\frac{1}{\omega_Nr^{N-1}}
                         \int\limits_{|x|=r}\left(1+u(x)\right)d\sigma_x$\\
                         &$\di=\;\Lambda \phi(r)\left(1+\bar u(r)\right)$\\
                         &$\di \leq\;\Lambda \phi(r)\left(1+U(r)\right),$\\
                         \end{tabular}$$
                         for all $\,r\geq 0$. It follows that
                         $$\di \left(r^{N-1}e^r\bar u'\right)'\leq\;\Lambda
                         e^rr^{N-1}\phi(r)\left(1+U(r)\right),
                         \qquad\mbox{for all}\;\,r\geq 0.$$
                         So, for all $\,r\geq r_0>0\,$,
                         $$\di\bar u(r)\leq\bar u(r_0)+\Lambda
                         \int_{r_0}^re^{-t}t^{1-N}\int_0^te^ss^{N-1}\phi(s)(1+U(s))dsdt.$$
                         The monotonicity of $U$ implies
                         \begin{equation}\label{azece}
                         \di\bar u(r)\leq \bar
                         u(r_0)+\Lambda(1+U(r))\int_{r_0}^r
                         e^{-t}t^{1-N}\int_0^te^ss^{N-1}\phi(s)dsdt,
                         \end{equation}
                         for all $r\geq r_0\geq 0.$
                          By \eq{phi} we can choose $r_0\geq 1$ such that
                         \begin{equation}\label{aunspe}
                         \di\int_{r_0}^{\infty}e^{-t}t^{1-N}\int_0^te^ss^{N-1}\phi(s)dsdt
                         <\frac{1}{2\Lambda}.
                         \end{equation}
                         Thus \eq{azece} and \eq{aunspe} yield
                         \begin{equation}\label{ff}
                         \di \bar u(r)\leq \bar u(r_0)+\frac{1}{2}(1+U(r)),
                         \qquad\mbox{for all}\;\;r\geq r_0.
                         \end{equation}
                         By the definition of $U$ and
                         $\di\lim_{r\ri\infty}\bar u(r)=\infty,$
                         we find $r_1\geq r_0$ such that
                         \begin{equation}\label{adoispe}
                         \di U(r)=\max_{r_0\leq t\leq r}\bar
                         u(r),\qquad\mbox{for all}\;\;r\geq r_1.
                         \end{equation}
                         Considering now \eq{ff} and \eq{adoispe} we obtain
                         $$\di U(r)\leq \bar
                         u(r_0)+\frac{1}{2}(1+U(r)),\qquad\mbox{for all}\;\;r\geq
                         r_1. $$
                         Hence
                         $$\di U(r)\leq 2\bar u(r_0)+1,\qquad\mbox{for
                         all}\;\;r\geq r_1.$$
                         This means that $U$ is bounded, so $u$ is also
                         bounded, a contradiction. It follows that \eq{unu} has
                         no
                         positive entire large solutions.

			 \medskip
			 \noindent{\sc Sufficient condition.}
                         We need the following auxiliary comparison result.

                         \begin{lemma} Assume that \eq{doi} and \eq{apatru}
                         hold.
                         Then the equations
                         \neweq{cinci}
                         \Delta v+|\nabla v|=\phi(|x|)f(v) \qquad
                         \Delta w+|\nabla w|=\psi(|x|)f(w)\\
                         \endeq
                         have positive entire large solution such that
                         \neweq{sase}
                         v\leq w  \qquad in\;\; \RR^N.\\
                         \endeq
                         \end{lemma}

                         {\bf Proof}. Radial solutions of \eq{cinci} satisfy
                         $$v''+\frac{N-1}{r}v'+|v'|=\phi(r)f(v)$$ and
                         $$w''+\frac{N-1}{r}w'+|w'|=\psi(r)f(w).$$
                         Assuming that $v'$ and $w'$ are non-negative, we
                         deduce
                         $$\displaystyle
                         \left(e^rr^{N-1}v'\right)'=e^rr^{N-1}\phi(r)f(v)$$
                         and
                         $$\di\left(e^rr^{N-1}w'\right)'=e^rr^{N-1}\psi(r)f(w).$$
                         Thus any positive solutions $v$ and $w$ of the
                         integral equations
                         \neweq{sapte}
                         v(r)=1+\int\limits_{0}^{r}e^{-t}t^{1-N}
                         \int\limits_{0}^{t}e^ss^{N-1}\phi(s)f(v(s))dsdt,\qquad
                         r\geq 0,\\
                         \endeq
                         \neweq{opt}
                         w(r)=b+\int\limits_{0}^{r}e^{-t}t^{1-N}
                         \int\limits_{0}^{t}e^ss^{N-1}\psi(s)f(w(s))dsdt,\qquad
                         r\geq 0,\\
                         \endeq
                         provide a solution of \eq{cinci}, for any $b>0$.
                         Since $w\geq b$, it follows that $f(w)\geq f(b)>0$
                         which yields
                         $$w(r)\geq b+f(b)\int\limits_{0}^{r}e^{-t}t^{1-N}
                         \int\limits_{0}^{t}e^ss^{N-1}\psi(s)dsdt,\qquad r\geq
                         0.\\$$
                         By \eq{apatru}, the right hand side of this
                         inequality goes to $+\infty$ as
                         $r\rightarrow \infty$. Thus $w(r)\ri\infty$ as
                         $r\ri\infty.$
                         With a similar argument we find $v(r)\ri\infty$ as
                         $r\ri\infty.$

                           Let $b>1$ be fixed.
                         We first show that \eq{opt} has a positive solution.
                         Similarly,
                         \eq{sapte} has a positive solution.
                         \medskip

                         Let $\{w_k\}$ be the sequence defined by $w_1=b$ and
                         \neweq{noua}
                         w_{k+1}(r)=b+\int\limits_{0}^{r}e^{-t}t^{1-N}
                         \int\limits_{0}^{t}e^ss^{N-1}\psi(s)f(w_k(s))dsdt,\qquad
                         k\geq 1.\\
                         \endeq

                         We remark that $\{w_k\}$ is a non-decreasing
                         sequence.
                         To get the convergence of $\{w_k\}$ we will show
                         that $\{w_k\}$ is bounded from above on bounded
                         subsets.
                         To this aim, we fix $R>0$ and we prove that
			  \neweq{zece}
                         w_k(r)\leq be^{Mr},\qquad\mbox{for any } 0\leq r\leq R,\;\mbox{ and for all
                         } k\geq 1,
                         \endeq
                         where $\di M\equiv\Lambda_N
                         \max_{t\in[0,R]}\,t\psi(t).$

                         We achieve \eq{zece} by induction. We first notice
                         that \eq{zece}
                         is true for $k=1$. Furthermore, the assumption $(f2)$
                         and the fact that $w_k\geq 1$
                         lead us to $f(w_k)\leq\Lambda w_k$, for all $k\geq
                         1$.
                         So, by \eq{noua},
                         $$w_{k+1}(r)\leq
                         b+\Lambda\int\limits_{0}^{r}e^{-t}t^{1-N}
                         \int\limits_{0}^{t}e^ss^{N-1}\psi(s)w_k(s)dsdt, \qquad
                         r\geq 0.$$\\
                         Using now \eq{doispe} (for
                         $g(t)=\psi(t)w_k(t)$) we deduce
                         $$w_{k+1}(r)\leq b+\Lambda
                         _N\int\limits_{0}^{r}t\psi(t)w_k(t)dt,
                         \qquad\forall\;r\in[0,R].$$
			 The induction hypothesis yields
                         $$w_{k+1}(r)\leq
                         b+bM\int\limits_{0}^{r}e^{Mt}dt=be^{Mr},\qquad\forall\;r\in[0,R].$$
			 Hence, by induction, the sequence
                         $\{w_k\}$ is
                         bounded in $[0,R]$, for any $R>0$.
                         It follows that $\di w(r)=\lim_{k\ri\infty}w_k(r)$ is
                         a positive solution of
                          \eq{opt}. In a similar way we conclude that
                         \eq{sapte} has a
                         positive solution on $[0,\infty)$.

                         The next step is to show that the constant $b$ may be
                         chosen sufficiently large
                         so that \eq{sase} holds. More exactly, if
                         \neweq{treispe}
                         b>1+K\Lambda_N\int\limits_{0}^{\infty}sh(s)\Psi(s)ds,
                         \endeq
                         where
                         $K=\exp\left(\Lambda_N\int\limits_{0}^{\infty}th(t)dt\right),$
                         then \eq{sase} occurs.

                         We first prove that the solution $v$ of \eq{sapte}
                         satisfies
                         \neweq{paispe}
                         v(r)\leq K\Psi(r),\qquad \forall\;r\geq 0.\\
                         \endeq
                         Since $v\geq 1$, from $(f2)$ we have $f(v)\leq \Lambda
                         v$.
                         We use this fact in \eq{sapte} and then we apply the
                         estimate
                         \eq{doispe} for $g=\phi.$ It follows that
                         \neweq{cincispe}
                         v(r)\leq
                         1+\Lambda_N\int\limits_{0}^{r}s\phi(s)v(s)ds,\qquad\forall\;r\geq
                         0.
                         \endeq
                         By Gronwall's inequality we obtain
                         $$v(r)\leq
                         \exp\left(\Lambda_N\int\limits_{0}^{r}s\phi(s)ds\right),
                         \qquad\forall\;r\geq 0,$$
                         and, by \eq{cincispe},
                         $$v(r)\leq 1+\Lambda_N\int\limits_{0}^{r}s\phi(s)\exp
                         \left(\Lambda_N\int\limits_{0}^{s}t\phi(t)dt\right)ds,
                         \qquad\forall\;r\geq 0.$$
                         Hence
                         $$v(r)\leq 1+\int\limits_{0}^{r}\left(\exp
                         \left(\Lambda_N\int\limits_{0}^{s}t\phi(t)dt\right)\right)'ds,
                         \qquad\forall\;r\geq 0,$$
                         that is
                         \neweq{58}
                         \di
                         v(r)\leq\exp\left(\Lambda_N\int\limits_{0}^{r}t\phi(t)dt\right),
                         \qquad\forall\;r\geq 0.
                         \endeq
                         Inserting $\phi=h+\psi$ in \eq{58} we have
                         $$\di v(r)\leq
                         e^{\Lambda_N\int\limits_{0}^{r}th(t)dt}\Psi(r)\leq
                         K\Psi(r),\qquad\forall\;r\geq 0,$$
                         so \eq{paispe} follows.

                         Since $b>1$ it follows that $v(0)<w(0).$ Then there
                         exists $R>0$
                         such that $v(r)<w(r),$ for any $0\leq r\leq R$. Set
                         $$R_{\infty}=\sup\{\ R>0 \,|\,v(r)<w(r),\;\,
                         \forall\,r\in[\,0,R] \,\}.$$
                         In order to conclude our proof, it remains to
                         show that $R_{\infty}=\infty$.
                         Suppose the contrary. Since
                         $v\leq w$ on $[\,0,R_{\infty}]$ and $\phi=h+\psi,$
                         from \eq{sapte} we deduce
                         $$\begin{tabular}{ll}
                         $\di v(R_{\infty})$ & $=\di
                         1+\int\limits_{0}^{\,R_{\infty}}e^{-t}t^{1-N}
                         \int\limits_{0}^{t}e^ss^{N-1}h(s)f(v(s))dsdt$\\
                         & $\di\qquad+\int\limits_{0}^{R_{\infty}}e^{-t}t^{1-N}
                         \int\limits_{0}^{t}e^ss^{N-1}\psi(s)f(v(s))dsdt.$\\
                         \end{tabular}$$
                         So, by \eq{doispe},
                         $$\di v(R_{\infty})
                         \leq
                         1+\frac{1}{N-2}\int\limits_{0}^{R_{\infty}}th(t)f(v(t))dt
                         +\int\limits_{0}^{R_{\infty}}e^{-t}t^{1-N}
                         \int\limits_{0}^{t}e^ss^{N-1}\psi(s)f(w(s))dsdt.$$
                         Taking into account that $v\geq 1$ and  the
                         assumption $(f2),$
                         it follows that
                         $$ v(R_{\infty})
                         \leq
                         1+K\Lambda_N\int\limits_{0}^{R_{\infty}}th(t)\Psi(t)dt
                         +\int\limits_{0}^{R_{\infty}}e^{-t}t^{1-N}
                         \int\limits_{0}^{t}e^ss^{N-1}\psi(s)f(w(s))dsdt.$$
                         Now, using \eq{treispe} we obtain
                         $$\di
                         v(R_{\infty})<b+\int\limits_{0}^{R_{\infty}}e^{-t}t^{1-N}
                         \int\limits_{0}^{t}e^ss^{N-1}\psi(s)f(w(s))dsdt=w(R_{\infty}).$$
                         Hence $v(R_{\infty})<w(R_{\infty}).$ Therefore,
                         there exists $R>R_{\infty}$ such that $v<w$ on
                         $[\,0,R]$, which
                         contradicts the maximality of $R_{\infty}$.
                         This contradiction shows that inequality \eq{sase}
                         holds and the
                         proof of Lemma 1 is now complete.
                         \endproof
                         \medskip

                         {\bf Proof of Theorem 2 completed.}
                         Suppose that \eq{apatru} holds. For all $k\geq 1$ we
                         consider the problem
                         \neweq{anouaspe}
                         \left\{\begin{tabular}{ll}
                         $\di\Delta u_k+|\nabla u_k|=p(x)f(u_k)$&$\di \mbox{ in
                         }\,\, B(0,k),$\\
                         $\di u_k(x)=w(k)$&$\di\mbox{ on }\,\,\partial
                         B(0,k).$\\
                         \end{tabular}\right.
                         \end{equation}
                         Then $v$ and $w$ defined by
                         \eq{sapte}
                         and \eq{opt} are positive sub and super--solutions of
                        \eq{anouaspe}.
                         So this problem has at least a positive solution
                         $u_k$ and
                         $$\di v(|x|)\leq u_k(x)\leq w(|x|) \qquad\mbox{ in
                         }\;B(0,k), \mbox{ for all }\;k\geq 1.$$
                         By Theorem 14.3 in \cite{gt}, the
                         sequence
                         $\{\nabla u_k\}$ is bounded on every compact set in
                         $\RR^N$.
                         Hence the sequence $\{u_k\}$ is bounded and
                         equicontinuous on compact subsets of
                         $\RR^N.$
                         So, by the Arzela-Ascoli Theorem, the sequence
                         $\{u_k\}$
                         has a uniform convergent subsequence, $\{u_k^1\}$ on
                         the ball
                         $B(0,1).$
                         Let $u^1=\lim_{k\rightarrow\infty}u_k^1$. Then
                         $\{f(u_k^1)\}$ converges
                         uniformly to $f(u^1)$ on  $B(0,1)$ and, by
                         \eq{anouaspe}, the sequence
                          $\{\Delta u_k^1+|\nabla u_k^1|\}$ converges uniformly
                         to
                         $pf(u^1).$
                         Since the sum of Laplacian and Gradient is a closed
                         operator, we deduce
                         that $u^1$ satisfies \eq{unu} on $B(0,1).$

                         Now, the sequence $\{u_k^1\}$ is bounded and
                         equicontinuous on the ball $B(0,2)$,  so it  has a
                         convergent subsequence $\{u_k^2\}.$
                         Let $u^2=\lim \limits _{k\rightarrow\infty}u_k^2\,$ on
                         $\,B(0,2)\,$
                         and $u^2$ satisfies \eq{unu} on $B(0,2).$ Proceeding
                         in the same way, we
                         construct a sequence $\{u^n\}$ so that $u^n$ satisfies
                         \eq{unu} on
                         $B(0,n)$
                         and $u^{n+1}=u^n$ on $B(0,n)$ for all $n$. Moreover,
                         the
                         sequence $\{u^n\}$
                         converges  in $L^\infty_{\rm loc}(\RR^N)$ to the
                         function $u$
                         defined by
                          $$u(x)=u^m(x), \qquad\mbox{for }\,\,x\in B(0,m).$$
                          Since $v\leq u^n\leq w$ on $B(0,n)$ it follows that
                         $v\leq u\leq
                         w$ on $\RR^N,$ and $\,u\,$ satisfies \eq{unu}. From
                         $v\leq u$ we
                         deduce that $u$ is a positive entire large solution of
                         \eq{unu}.
                         This completes the proof.
                         \endproof

                         \end{document}